\documentclass[12pt,a4paper]{article}
\usepackage{}
\usepackage{indentfirst}
\usepackage{amsfonts}
\usepackage{amssymb}
\usepackage{mathrsfs}
\usepackage{amsmath}
\usepackage{amsthm}
\usepackage{enumerate}
\usepackage{cite}
\usepackage[all]{xy}
\allowdisplaybreaks
\newtheorem{defn}{Definition}[section]
\newtheorem{thm}[defn]{Theorem}
\newtheorem{lem}[defn]{Lemma}
\newtheorem{prop}[defn]{Proposition}
\newtheorem{cor}[defn]{Corollary}
\newtheorem{eg}[defn]{Example}
\newtheorem{re}[defn]{Remark}
\newcommand{\bdefn}{\begin{defn}}
\newcommand{\edefn}{\end{defn}}
\newcommand{\bthm}{\begin{thm}}
\newcommand{\ethm}{\end{thm}}
\newcommand{\blem}{\begin{lem}}
\newcommand{\elem}{\end{lem}}
\newcommand{\bprop}{\begin{prop}}
\newcommand{\eprop}{\end{prop}}
\newcommand{\bcor}{\begin{cor}}
\newcommand{\ecor}{\end{cor}}
\newcommand{\beg}{\begin{eg}}
\newcommand{\eeg}{\end{eg}}
\newcommand{\bre}{\begin{re}}
\newcommand{\ere}{\end{re}}
\newcommand{\bpf}{\begin{proof}}
\newcommand{\epf}{\end{proof}}

\newcommand{\Id}{{\rm Id}}

\newcommand{\Ker}{\rm Ker}

\newcommand{\K}{\mathbb{K}}

\newcommand{\supercite}[1]{\textsuperscript{\cite{#1}}}
\newcommand{\benu}{\begin{enumerate}}
\newcommand{\eenu}{\end{enumerate}}
\newcommand{\bc}{\begin{center}}
\newcommand{\ec}{\end{center}}
\newcommand{\bea}{\begin{eqnarray}}
\newcommand{\eea}{\end{eqnarray}}
\newcommand{\Bea}{\begin{eqnarray*}}
\newcommand{\Eea}{\end{eqnarray*}}
\newcommand{\beq}{\begin{equation}}
\newcommand{\eeq}{\end{equation}}
\newcommand{\Beq}{\begin{equation*}}
\newcommand{\Eeq}{\end{equation*}}
\newcommand{\bspl}{\begin{split}}
\newcommand{\espl}{\end{split}}
\newcommand\relphantom[1]{\mathrel{\phantom{#1}}}

\numberwithin{equation}{section}

\begin{document}
\title{\bf  On universal $\alpha$-central extensions of Hom-preLie algebras}
\author{\normalsize \bf Bing Sun$^1$, Liangyun Chen$^2$, Xin Zhou$^{2,3}$}
\date{{\small {$^1$ ．College of Mathematics, Changchun Normal
 University,\\
Changchun 130032, China}\\
\small {$^2$ School of
Mathematics and Statistics, Northeast Normal
 University,\\Changchun 130024, China} \\
\small {$^3$ School of Mathematical and Statistics, Yili Normal University,\\ Yining 835000, China} }}
\maketitle
\begin{abstract}
We introduce the notion of Hom-co-represention and low-dimensional chain complex. We study universal central extensions of Hom-preLie algebras and generlize some classical results. As the same time, we introduce $\alpha$-central extensions, universal $\alpha$-central extensions and $\alpha$-perfect Hom-preLie algebras. We construct universal ($\alpha$)-central extensions of Hom-preLie algebras.
\bigskip

\end{abstract}

\footnote[0]{Corresponding author(L. Chen): chenly640@nenu.edu.cn.}
\footnote[0]{Supported by  NNSF of China (Nos. 11171055 and 11471090). }

\section{Introduction}

A Hom-preLie algebra, as a kind of Hom-Lie admissible algebra, was introduced by
Makhlouf-Silvestrov \cite{Makhlouf1}.
Specifically,
for a vector space $L$ over a field $\K$  equipped with a bilinear map $\mu :L\times L\rightarrow L$ and a linear map $\alpha:L\rightarrow L$, we say that the triple
$(L,\mu,\alpha)$ is a \textit{Hom-preLie algebra} if
\begin{align*}
\alpha(x)(yz)-(xy)\alpha(z)=\alpha(y)(xz)-(yx)\alpha(z),
\end{align*}
for all $x,y,z\in L$. If the elements of $L$ also satisfy the following equation
\begin{align*}
(xy)\alpha(z)=(xz)\alpha(y).
\end{align*}
Then we call $(L,\mu,\alpha)$ is a Hom-Novikov algebra.
Clearly, a Hom-Novikov algebra is a Hom-preLie algebra. Moreover, Hom-preLie algebras generalizes the notation of pre-Lie algebras ($\alpha=\Id_{L}$), which has been extensively studied in the construction and classification of Hom-Novikov algebras\supercite{Yau}.





In recent year, the universal central extension
of a perfect Leibniz algebra was studied in several articles\supercite{Casas1,Casas2,Casas3,Casas4,Gnedbaye1,Gnedbaye2,Kurdiani}. In \cite{Casas5,Casas6,Khmaladze},
authors  study universal ($\alpha$)-central extension.
The purpose of this paper is to study universal $\alpha$-central extensions of Hom-preLie algebras by the way of \cite{Casas5}

The organization of this paper is as follows: after this introduction, in Section 2 we recall the
background material on Hom-preLie algebras. In section 3, we define Hom-co-representations and  low-dimensional chain complex, which derive a low-dimensional homology $\K$-vector space of a Hom-preLie algebra. In section 4, we extend classical results and present a counterexample showing that the composition of two central extensions is not a central extensionon. The fact lead us to define $\alpha$-central extensions. Moreover, We construct a right exact covariant functor $\frak{uce}$ of a Hom-preLie algebra which acts on a perfect Hom-preLie algebra $L$, obtain the universal central extension of the perfect Hom-preLie algebra $\frak{uce}_\alpha(L)=\frac{\alpha_L(L)\otimes \alpha_L(L)}{I_L}$, where $I_L=\langle\alpha_L(x_1)\otimes x_2x_3- x_1x_2\otimes \alpha_L(x_3)-\alpha_L(x_2)\otimes x_1x_3+ x_2x_1\otimes \alpha_L(x_3)\rangle.$ And the kernel of the univeral central extension is exactly the second homology with trivial coefficients of the Hom-preLie algebra. We also introduce the notion of $\alpha$-perfect Hom-preLie algebra and then construct the universal $\alpha$-central extension of $\alpha$-perfect Hom-preLie algebra.

\section{Preliminaries on Hom-preLie algebras}

Throughout this paper $\K$ denotes an arbitrary field.
\bdefn\rm
A homomorphism of Hom-preLie algebras $f:(L,\mu,\alpha_L) \to (L',\mu',\alpha_{L'})$ is a $\mathbb{K}$-linear map $f : L \to L'$ such that
\begin{enumerate}
\item[a)] $f(\mu(x,y)) =\mu'(f(x),f(y)),$
\item [b)] $f \circ \alpha_L(x) = \alpha_{L'} \circ f(x),$
\end{enumerate}
for all $x, y \in L$.

The Hom-preLie algebras $(L,\mu,\alpha_L)$ and $(L',\mu',\alpha_{L'})$ are isomorphic if there is a  Hom-preLie algebras homomorphism $f: (L,\mu,\alpha_L) \to (L',\mu',\alpha_{L'})$ such that $f:L \to L'$ is bijective.
\edefn

\bdefn\rm
Let $(L,\alpha_L)$ be a Hom-preLie algebra. A  Hom-preLie subalgebra $(H, \alpha_H)$ of $(L,\alpha_L)$ is a linear subspace $H$ of $L$, which is closed for the bracket $\mu$ and invariant by $\alpha_L$, that is,
\begin{enumerate}
\item [a)] $xy \in H,$ for all $x, y \in H$,
\item [b)] $\alpha_L(x) \in H$, for all $x \in H$,
\end{enumerate}
where $\mu_H$ and $\alpha_H$ for $H$ are the restrictions of $\mu_L$ and $\alpha_L$.

A  Hom-preLie subalgebra $(H, \alpha_H)$ of $(L,\alpha_L)$ is said to be a  Hom-ideal if $xy,yx \in H$ for all $x \in H, y \in L$.
\edefn
If $(H, \alpha_H)$ is a  Hom-ideal of $(L,\alpha_L)$, then  $(L/H, \overline{\alpha_L})$ naturally inherits a structure of Hom-preLie algebra, where $\overline{\alpha_L}: L/H\rightarrow L/H $ is a linear map, which is said to be the quotient Hom-preLie algebra.
\bdefn
The annihilator of a  Hom-preLie algebra  $(L,\alpha_L)$ is the $\mathbb{K}$-vector subspace  $$Z(L) = \{ x \in L \mid xy=0=yx , for\ all\ y \in L\}.$$
\edefn

\bdefn\rm
 A short exact sequence of Hom-preLie algebras  $(K) : 0 \to (M, \alpha_M) \stackrel{i} \to (K,\alpha_K) \stackrel{\pi} \to (L, \alpha_L) \to 0$ is said to be split if
 there exists a Hom-preLie algebra homomorphism $\sigma : (L, \alpha_L) \to (K, \alpha_K)$ such that $\pi \circ \sigma = {\Id}_L$.
\edefn

Let $(M, \alpha_M)$ and $(L,\alpha_L)$ be Hom-preLie algebras with a Hom-action of $(L,\alpha_L)$ over $(M, \alpha_M)$. Then we can construct the sequence
 $$0 \to (M,\alpha_M) \stackrel{i}\to (M\rtimes L,\tilde{\alpha}) \stackrel{\pi}\to (L,\alpha_L) \to 0,$$
where $i:M \to M\rtimes L,i(m)=(m,0)$ and $\pi:M\rtimes L \to L, \pi(m,l)=l$. Moreover, this sequence splits by $\sigma: L \to M\rtimes L, \sigma(l)=(0,l)$, that is, $\sigma$ satisfies $\pi \circ \sigma={\Id}_L$ and $\tilde{\alpha}\circ \sigma=\sigma \circ \alpha_L$.

\bdefn\rm
Let $(M, \alpha_M)$ be a Hom-representation of a Hom-preLie algebra
$(L, \alpha_L)$. A derivation of  $(L, \alpha_L)$ over $(M, \alpha_M)$ is a $\K$-linear map $d : L \to M$
satisfying:
\begin{enumerate}
 \item[a)] $d(l_1l_2)=\alpha_L(l_1)\centerdot d(l_2)+d(l_1)\centerdot \alpha_L(l_2),$
 \item [b)] $d\circ \alpha_L=\alpha_M\circ d,$
 \end{enumerate}
 for all $l_1,l_2\in L$.
\edefn

\section{Homology of Hom-preLie algebra}

\bdefn\rm
Let $(L,\alpha_L)$ be a Hom-preLie algebra.
A Hom-co-representation of $(L,\alpha_L)$ is a $\K$-vector space $M$ together with two bilinear maps, $\lambda: L \otimes M \to M,$ $\lambda(l\otimes m)=l\centerdot m$, and $\rho: M \otimes L \to M,$ $\rho(m\otimes l)=m\centerdot l$, and a $\K$-linear map $\alpha_M: M \to M$ satisfying the following identities:
\begin{enumerate}
 \item[a)] $\alpha_M(m)\centerdot (xy) = (m\centerdot x) \centerdot \alpha_L(y) ,$
      \item [b)] $(xy)\centerdot \alpha_M(m)=\alpha_L(x)\centerdot (y\centerdot m),$
 \item [c)] $\alpha_L(x)\centerdot
     (m\centerdot y)=(x\centerdot m)\centerdot \alpha_L(y),$
 \item [d)] $\alpha_M(x\centerdot m)= \alpha_L(x) \centerdot \alpha_M(m),$
 \item [e)] $\alpha_M(m\centerdot x)=\alpha_M(m) \centerdot  \alpha_L(x),$
 \end{enumerate}
 for all $x, y \in L$ and $m, m' \in M$.
\edefn

\beg\rm\label{eg212}
Every Hom-preLie algebra $(L,\alpha_L)$ has a Hom-co-representation structure on itself given by the actions
$$x\centerdot m = mx; \ \ m\centerdot x =mx,$$
where $x\in L$ and m is an element of the underlying $\K$-vector space to $L$.
\eeg

Let $(L,\alpha_L)$ be a Hom-preLie algebra and $(M, \alpha_M)$ be a Hom-co-representation of $(L,\alpha_L)$. Denote $CL_n^\alpha (L, M) := M \otimes L^{\otimes n}, n\geq 0\ (CL_{-1}^\alpha (L, M) =0).$ We
define the $\K$-linear maps
$$d_n : CL_n^\alpha (L, M)\to CL_{n-1}^\alpha (L, M)(n=0,1,2,3)$$
by
\begin{align*}
d_0(m)=&0;\\
d_1(m\otimes x_1)=&m\centerdot x_1-x_1\centerdot m;\\
d_2(m\otimes x_1\otimes x_2)=&m\centerdot x_1\otimes \alpha_L(x_2)+x_2\centerdot m\otimes \alpha_L(x_1)-\alpha_M(m)\otimes x_1x_2;\\
d_3(m\otimes x_1\otimes x_2\otimes x_3)=&\alpha_M(m)\otimes \alpha_L(x_1)\otimes x_2x_3-\alpha_M(m)\otimes x_1x_2\otimes \alpha_L(x_3)\\
&-\alpha_M(m)\otimes \alpha_L(x_2)\otimes x_1x_3+\alpha_M(m)\otimes x_2x_1\otimes \alpha_L(x_3)\\
&+m\centerdot x_1\otimes \alpha_L(x_2)\otimes \alpha_L(x_3)-x_3\centerdot m\otimes \alpha_L(x_1)\otimes \alpha_L(x_2)\\
&-m\centerdot x_2\otimes \alpha_L(x_1)\otimes \alpha_L(x_3)+x_3\centerdot m\otimes \alpha_L(x_2)\otimes \alpha_L(x_1).
\end{align*}

\bprop\label{prop212}
With the above notations, we have $d_{n-1}\circ d_n=0\ (n=2,3)$.
\eprop
\bpf
Since $(M, \alpha_M)$ is a Hom-co-representation of $(L,\alpha_L)$, we have
\begin{align*}
d_{1}\circ d_2(m\otimes x_1\otimes x_2)
=&(m\centerdot x_1)\centerdot \alpha_L(x_2)-\alpha_L(x_2)\centerdot (m\centerdot x_1)+(x_2\centerdot m)\centerdot \alpha_L(x_1)\\
&-\alpha_L(x_1)\centerdot (x_2\centerdot m)-\alpha_M(m)\centerdot (x_1x_2)+(x_1x_2)\centerdot \alpha_M(m)\\
=&0
\end{align*}
and
\begin{align*}
&d_{2}\circ d_3(m\otimes x_1\otimes x_2\otimes x_3)\\
=&\alpha_M(m)\centerdot \alpha_L(x_1)\otimes \alpha_L(x_2x_3)+(x_2x_3)\centerdot \alpha_M(m)\otimes \alpha^2_L(x_1)-\alpha^2_M(m)\otimes \alpha_L(x_1)(x_2x_3)\\
&-\alpha_M(m)\centerdot (x_1x_2)\otimes \alpha^2_L(x_3)-\alpha_L(x_3)\centerdot \alpha_M(m)\otimes \alpha_L(x_1x_2)+\alpha^2_M(m)\otimes(x_1x_2)\alpha_L(x_3)\\
&-\alpha_M(m)\centerdot \alpha_L(x_2)\otimes \alpha_L(x_1x_3)-(x_1x_3)\centerdot \alpha_M(m)\otimes  \alpha^2_L(x_2)+\alpha^2_M(m)\otimes \alpha_L(x_2)(x_1x_3)\\
&+\alpha_M(m)\centerdot (x_2x_1)\otimes \alpha^2_L(x_3)+\alpha_L(x_3)\centerdot \alpha_M(m)\otimes \alpha_L(x_2x_1)-\alpha^2_M(m)\otimes (x_2x_1)\alpha_L(x_3)\\
&+(m\centerdot x_1)\centerdot \alpha_L(x_2)\otimes \alpha^2_L(x_3)\!+\!\alpha_L(x_3)\centerdot (m\centerdot x_1)\otimes \alpha^2_L(x_2)\!-\!\alpha_M(m\centerdot x_1)\otimes\alpha_L(x_2)\alpha_L(x_3)\\
&-(x_3\centerdot m)\centerdot \alpha_L(x_1)\otimes \alpha^2_L(x_2)\!-\!\alpha_L(x_2)\centerdot (x_3\centerdot m)\otimes \alpha^2_L(x_1)\!+\! \alpha_M(x_3\centerdot m)\otimes  \alpha_L(x_1)\alpha_L(x_2)\\
&-( m\centerdot x_2)\centerdot \alpha_L(x_1)\otimes \alpha^2_L(x_3)\!-\!\alpha_L(x_3)\centerdot (m\centerdot x_2)\otimes \alpha^2_L(x_1)\!+\!\alpha_M(m\centerdot x_2)\otimes \alpha_L(x_1)\alpha_L(x_3)\\
&+(x_3\centerdot m)\centerdot \alpha_L(x_2)\otimes \alpha^2_L(x_1)\!+\!\alpha_L(x_1)\centerdot (x_3\centerdot m)\otimes \alpha^2_L(x_2)\!- \! \alpha_M(x_3\centerdot m)\otimes  \alpha_L(x_2)\alpha_L(x_1)\\
=&0.
\end{align*}
Hence, we prove this proposition.
\epf
 By the Proposition \ref{prop212}, we can define the chain complex $(CL_\star^{\alpha},d_\star)$. The homology of the chain complex $(CL_\star^{\alpha},d_\star)$ is called homology of the Hom-preLie algebra $(L, \alpha_L)$ with coefficients in the Hom-co-representation $(M, \alpha_M)$ and we will denote it by
$$HL_\star^\alpha(L, M):= H_\star(CL_\star^\alpha(L, M), d_\star).$$
Now we are going to compute low dimensional homologies. So, for n = 1, a
direct checking shows that
$$HL_0^\alpha(L, M)=\frac{M}{M_L},$$
where $M_L=\{m\centerdot l-l\centerdot m: m\in M,l\in L\}.$

  If $M$ is a trivial Hom-co-representation and $n=2$, that is, $m\centerdot l=l\centerdot m=0$, then
  $$HL_1^\alpha(L, M)=\frac{M\otimes L}{\alpha_M(M)\otimes (LL)}.$$

\bprop
Let $(L, \alpha_L)$ be a Hom-preLie algebra, which is considered as a Hom-co-representation of itself as in Example \ref{eg212} and $\K$ as a trivial
Hom-co-representation of  $(L, \alpha_L)$. Then
$$HL_n^\alpha(L, L)\cong HL_{n+1}^\alpha(L, \K),(n=0,1).$$
\eprop
\bpf
Let $-{\Id}$ be a minus identity transformation.
Obviously $-{\Id} : CL_n^\alpha(L, L)\to CL_{n+1}^\alpha(L, \K), (n = 0,1)$ defines a chain isomorphism, hence the isomorphism in the homologies.
\epf

\section{Universal ($\alpha$)-central extensions}

\bdefn\rm
 A short exact sequence of Hom-preLie algebras  $(K) : 0 \to (M, \alpha_M) \stackrel{i} \to (K,\alpha_K) \stackrel{\pi} \to (L, \alpha_L) \to 0$ is said to be central if $MK= 0 =KM$, equivalently,  $M \subseteq Z(K)$.

 The sequence $(K)$ is said to be $\alpha$-central if $\alpha_M(M)K = 0=K \alpha_M(M)$, equivalently, $\alpha_M(M) \subseteq Z(K)$.
\edefn

\bre\label{re211}
Obviously, every  central extension is an $\alpha$-central extension. Note that in the case $\alpha_M = {\Id}_M$, both notions coincide.
\ere

\bdefn\rm
A central extension  $(K) : 0 \to (M, \alpha_M) \stackrel{i} \to (K,\alpha_K) \stackrel{\pi} \to (L, \alpha_L) \to 0$ is said to be universal if for every   central extension $(K') : 0 \to (M', \alpha_{M'}) \stackrel{i'} \to (K',\alpha_{K'}) \stackrel{\pi'} \to (L, \alpha_L) \to 0$  there exists a  unique homomorphism of Hom-preLie algebras  $h : (K,\alpha_K) \to (K',\alpha_{K'})$ such that $\pi'\circ h = \pi$.

A central extension  $(K) : 0 \to (M, \alpha_M) \stackrel{i} \to (K,\alpha_K) \stackrel{\pi} \to (L, \alpha_L) \to 0$ is said to be universal $\alpha$-central  if for every  $\alpha$-central extension $(K') : 0 \to (M', \alpha_{M'}) \stackrel{i'} \to (K',\alpha_{K'}) \stackrel{\pi'} \to (L, \alpha_L) \to 0$  there exists a unique homomorphism of  Hom-preLie algebras  $h : (K,\alpha_K) \to (K',\alpha_{K'})$ such that $\pi' \circ h = \pi$.
\edefn

\bre\rm\label{re212}
Obviously, every universal $\alpha$-central extension  is a  universal  central extension.  Note that in the case $\alpha_M = {\Id}_M$, both notions coincide.
\ere

\bdefn\rm
A Hom-preLie algebra  $(L, \alpha_L)$ is said to be perfect if  $L = LL$, where $LL$ denotes the linear span of the elements $xy$, with $x,y\in L.$
\edefn

\blem\label{lem212}
Let $\pi : (K,\alpha_K) \to (L, \alpha_L)$ be a surjective homomorphism  of  Hom-preLie algebras. If $(K,\alpha_K)$ is a perfect  Hom-preLie algebra, then $(L, \alpha_L)$ is a perfect Hom-preLie algebra as well.
\elem
\bpf
Since $(K,\alpha_K)$ is a perfect  Hom-preLie algebra and $\pi$ is a surjective homomorphism, we have
$$LL=\pi(K)\pi(K)=\pi(KK)=\pi(K)=L.$$
The proof is completed.
\epf

\blem\label{lem213}
Let $0 \to (M, \alpha_M) \stackrel{i} \to (K,\alpha_K) \stackrel{\pi} \to (L, \alpha_L) \to 0$ be an  $\alpha$-central extension and $(K,\alpha_K)$ a perfect Hom-preLie algebra. If there exists a homomorphism of Hom-preLie algebras  $f : (K,\alpha_K) \to (A, \alpha_A)$ such that $\tau \circ f = \pi$, where $0 \to (N, \alpha_N) \stackrel{j} \to (A,\alpha_A) \stackrel{\tau} \to (L, \alpha_L) \to 0$ is a central extension, then $f$ is unique.
\elem
\bpf
Let us assume that there exist two homomorphisms $f_1$ and $f_2$ such that
$\tau \circ f_1=\pi=\tau \circ f_2$, then $f_1-f_2\in {\Ker} \tau = N$, i.e., $f_1(k) = f_2(k) + n_k, n_k\in N$
and $N\subseteq Z(A)$.
Moreover, $f_1$ and $f_2$ coincide over $KK$. Indeed
\begin{align*}
f_1(k_1k_2)=&f_1(k_1)f_1(k_2)=(f_2(k_1)+n_{k_1})(f_2(k_2)+n_{k_2})\\
=&f_2(k_1)f_2(k_2)+f_2(k_1)n_{k_2}+n_{k_1}f_2(k_2)+n_{k_1}n_{k_2}=f_2(k_1)f_2(k_2)\\
=&f_2(k_1k_2).
\end{align*}
Since $(K, \alpha_K)$ is perfect, $f_1$ coincides with $f_2$ over $K$.
\epf

\blem\label{lem214}
If $0 \to (M, \alpha_M) \stackrel{i} \to (K,\alpha_K) \stackrel{\pi} \to (L, \alpha_L) \to 0$ is a universal central extension, then  $(K,\alpha_K)$ and $(L, \alpha_L)$ are perfect Hom-preLie algebras.
\elem
\bpf
Let us assume that $(K,\alpha_K)$ is not a perfect Hom-preLie algebra, then $KK \varsubsetneq K$. Consider $I$ the smallest Hom-ideal generated by $KK$, then $(K/I, \widetilde{\alpha})$, where $\widetilde{\alpha}$ is the induced  homomorphism, is an abelian Hom-preLie algebra, consequently,  is a trivial Hom-co-representation of $(L,\alpha_L)$. Let us consider the central extension $0 \to (K/I, \widetilde{\alpha}) \to (K/I \times L, \widetilde{\alpha}\times \alpha_L) \stackrel{pr}\to  (L, \alpha_L) \to 0$. Then the  homomorphisms of Hom-preLie algebras  $\varphi, \psi : (K,\alpha_K) \to (K/I \times L, \widetilde{\alpha}\times \alpha_L)$ given by $\varphi(k)=(k+I,\pi(k))$ and $\psi(k)=(0,\pi(k)), k \in K$, verify that  $pr\circ  \varphi = \pi = pr \circ  \psi$, so $0 \to (M, \alpha_M) \stackrel{i} \to (K,\alpha_K) \stackrel{\pi} \to (L, \alpha_L) \to 0$ cannot be a universal central extension. Lemma \ref{lem212} ends the proof.
\epf

\beg\rm
Consider the two-dimensional Hom-preLie algebra $(L, \alpha_L)$ with the
basis ${b_1, b_2}$, the multiplication given by $b_2b_1= b_2, b_2b_2 = b_1$ (and all other products 0) and the endomorphism $\alpha_L = 0.$

Let $(K, \alpha_K)$ be the three-dimensional  Hom-preLie algebra with the basis ${a_1, a_2, a_3}$,
multiplication given by $a_2a_2 = a_1, a_3a_2 = a_3, a_3a_3 = a_2 $ (and all other products 0) and the endomorphism $\alpha_K = 0.$
Obviously $(K, \alpha_K)$ is a perfect  Hom-preLie algebra since $K = KK$ and
$Z(K) = \langle a_1\rangle.$
The linear map $\pi : (K, 0)\to (L, 0)$ given by $\pi(a_1) = 0, \pi(a_2) = b_1, \pi(a_3) = b_2,$
is a central extension since $\pi$ is trivially surjective and is a homomorphism of
Hom-preLie algebras:
\begin{alignat*}{3}
\pi(a_2a_2)&=\pi(a_1)=0;\,\,\, &\,\,\, \pi(a_2)\pi(a_2)&=b_1b_1=0,\\
\pi(a_2a_3)&=\pi(0)=0;\,\,\, &\,\,\, \pi(a_2)\pi(a_3)&=b_1b_2=0,\\
\pi(a_3a_2)&=\pi(a_3)=b_2; & \pi(a_3)\pi(a_2)&=b_2b_1=b_2,\\
\pi(a_3a_3)&=\pi(a_2)=b_1; & \pi(a_3)\pi(a_3)&=b_2b_2=b_1.
\end{alignat*}
Obviously, $0\circ \pi = \pi \circ 0$ and $\Ker$$(\pi) = \langle a_1\rangle$, hence $\Ker$$(\pi) \subseteq Z(K)$.

Now consider the four-dimensional Hom-preLie algebra $(F, \alpha_F)$ with the basis
${e_1, e_2, e_3, e_4}$, the multiplication given by $e_2e_3 = e_1, e_3e_3 = e_2, e_4e_3 = e_4, e_4e_4 = e_3$ (and all other products 0) and the endomorphism $\alpha_F = 0$.

The linear map $\rho : (F, 0) \!\to\! (K, 0)$ given by $\rho(e_1) \!= \!0, \rho(e_2) \!= \!a_1, \rho(e_3) \!=\!
a_2, \rho(e_4) \!=\! a_3,$ is a central extension since $\rho$ is trivially surjective and is a homomorphism of Hom-preLie algebras:
\begin{alignat*}{4}
\rho(e_2e_3)&=\rho(e_1)=0;\,\,\,\, &\,\,\,\, \rho(e_2)\rho(e_3)&=a_1a_2=0,\\
\rho(e_2e_4)&=\rho(0)=0;\,\,\,\, &\,\,\,\, \rho(e_2)\rho(e_4)&=a_1a_3=0,\\
\rho(e_3e_2)&=\rho(0)=0;\,\,\,\, &\,\,\,\, \rho(e_3)\rho(e_2)&=a_2a_1=0,\\
\rho(e_3e_3)&=\rho(e_2)=a_1; & \rho(e_3)\rho(e_3)&=a_2a_2=a_1,\\
\rho(e_3e_4)&=\rho(0)=0; & \rho(e_3)\rho(e_4)&=a_2a_3=0,\\
\rho(a_4e_3)&=\rho(e_4)=a_3; & \rho(e_4)\rho(e_3)&=a_3a_2=a_3,\\
\rho(e_4e_4)&=\rho(e_3)=a_2; & \rho(a_4)\rho(a_4)&=a_3a_3=a_2.
\end{alignat*}
Obviously, $0\circ \rho = \rho \circ 0$ and $\Ker$$(\rho) = \langle e_1 \rangle$ and $Z(F) = \langle e_1 \rangle$, hence $\Ker$$(\rho) \subseteq Z(F)$.

The composition $\pi\circ \rho : (F, 0) \to (L, 0)$ is given by
\begin{align*}
\pi\circ \rho (e_1) &= \pi(0) = 0, \ \ \ \ \ \pi\circ \rho (e_2) = \pi(a_1) = 0,\\
\pi\circ \rho (e_3) &= \pi(a_2) = b_1, \ \ \ \  \pi\circ \rho(e_4) = \pi(a_3) = b_2.
 \end{align*}
Consequently,
$\pi\circ \rho : (F, 0) \to (L, 0)$
is a surjective homomorphism, but is not a central extension,
since $Z(F) = \langle e_1\rangle$ and $\Ker$$(\pi\circ \rho) =\langle e_1,e_2 \rangle$, i.e., $\Ker$$(\pi\circ \rho)\nsubseteq Z(F )$.
\eeg
\blem\label{lem215}
Let $0 \to (M, \alpha_M) \stackrel{i} \to (K,\alpha_K) \stackrel{\pi} \to (L, \alpha_L) \to 0$  and  $0 \to (N, \alpha_N) \stackrel{j} \to (F,\alpha_F) \stackrel{\rho} \to (K, \alpha_K) \to 0$ be  central extensions with   $(K, \alpha_K)$ a perfect Hom-preLie algebra. Then the composition extension $0 \to (P, \alpha_P) = {\rm Ker}\ (\pi \circ \rho)   \to (F,\alpha_F) \stackrel{\pi \circ \rho} \to (L, \alpha_L) \to 0$ is an $\alpha$-central extension.

Moreover, if $0 \to (M, \alpha_M) \stackrel{i} \to (K,\alpha_K) \stackrel{\pi} \to (L, \alpha_L) \to 0$ is a universal  $\alpha$-central extension, then $0 \to (N, \alpha_N) \stackrel{j} \to (F,\alpha_F) \stackrel{\rho} \to (K, \alpha_K) \to 0$ is split.
\elem

\bpf
 We must prove that  $\alpha_P(P)F=F \alpha_P(P)=0$.

Since $(K,\alpha_K)$ is a perfect  Hom-preLie algebra, every element $f \in F$ can be written as  $f = \displaystyle \sum_i \lambda_i f_{i_1}f_{i_2}  + n, n \in N, \lambda_i \in \mathbb{K}, f_{i_j} \in F, j=1,2$.
So, for all $p \in P, f \in F$ we have that
\begin{align*}
\alpha_P(p) f =&\displaystyle \sum_i \lambda_i \alpha_P(p)(f_{i_1}f_{i_2})+\alpha_P(p)n\\=& \displaystyle \sum_i \lambda_i \left( (pf_{i_1})\alpha_F (f_{i_2}) +\alpha_F(f_{i_1})(pf_{i_2})-(f_{i_1}p)\alpha_F(f_{i_2}) \right)+ \alpha_P(p)n = 0,
\end{align*}
since $\rho(pf_{i_j})=\rho(p)\rho(f_{i_j})= {\rm Ker}(\pi)K =0$.

In a similar way we can check that $f \alpha_P(p) =0$.

We now consider the second statement. If $0 \to (M, \alpha_M) \stackrel{i} \to (K,\alpha_K) \stackrel{\pi} \to (L, \alpha_L) \to 0$ is a universal $\alpha$-central extension, then by the first statement, $0 \to (P, \alpha_P) = {\rm Ker}\ (\pi \circ \rho)   \to (F,\alpha_F) \stackrel{\pi \circ \rho} \to (L, \alpha_L) \to 0$ is an  $\alpha$-central extension, then there exists a unique  homomorphism of Hom-Lie algebras  $\sigma : (K,\alpha_K) \to (F,\alpha_F)$ such that $\pi \circ \rho \circ \sigma = \pi$. On the other hand, $\pi \circ \rho \circ \sigma = \pi = \pi \circ {\Id}$ and $(K,\alpha_K)$ is perfect, then Lemma \ref{lem213} implies that $\rho \circ \sigma = {\Id}$.
\epf

\bthm\label{thm211}
a)  If a central extension  $0 \to (M, \alpha_M) \stackrel{i} \to (K,\alpha_K) \stackrel{\pi} \to (L, \alpha_L) \to 0$ is a universal $\alpha$-central extension, then  $(K,\alpha_K)$ is a perfect   Hom-preLie algebra and every central extension  of $(K,\alpha_K)$ is split.

    b) Let $0 \to (M, \alpha_M) \stackrel{i} \to (K,\alpha_K) \stackrel{\pi} \to (L, \alpha_L) \to 0$ be a central extension.
If $(K,\alpha_K)$ is a perfect  Hom-preLie algebra and every  central extension of $(K,\alpha_K)$ is split, then $0 \to (M, \alpha_M) \stackrel{i} \to (K,\alpha_K) \stackrel{\pi} \to (L, \alpha_L) \to 0$ is a universal central extension.

c) A Hom-preLie algebra $(L, \alpha_L)$ admits a universal central extension if and only if $(L, \alpha_L)$ is perfect.

d) The kernel of the universal central extension is canonically isomorphic to $HL_2^{\alpha}(L)$.
\ethm
\bpf
 \noindent  a) If  $0 \to (M, \alpha_M) \stackrel{i} \to (K,\alpha_K) \stackrel{\pi} \to (L, \alpha_L) \to 0$ is a universal $\alpha$-central extension, then it is a universal central extension by  Remark \ref{re212}, so $(K,\alpha_K)$ is a perfect Hom-preLie algebra by Lemma \ref{lem214} and every  central extension of  $(K,\alpha_K)$ is split by Lemma \ref{lem215}.

   b) Consider a central extension $0 \to (N, \alpha_N) \stackrel{j} \to (A,\alpha_A) \stackrel{\tau} \to (L, \alpha_L) \to 0$. Construct the pull-back extension   $0 \to (N, \alpha_N) \stackrel{\chi} \to (P,\alpha_P) \stackrel{\overline{\tau}} \to (K, \alpha_K) \to 0$, where $P=\{(a,k) \in A \times K \mid \tau(a)=\pi(k) \}$ and  $\alpha_P(a,k)=(\alpha_A(a),\alpha_K(k))$, which is central. Since every  central extension of $(K,\alpha_K)$ is split, there exists a homomorphism $\sigma : (K,\alpha_K) \to (P,\alpha_P)$ such that $\overline{\tau} \circ \sigma = {\Id}$. Then $\overline{\tau} \circ \sigma$ By the pull-back construction, there exists a homomorphism  $\overline{\pi} :  (P,\alpha_P)\to (A,\alpha_A)$ such that the following diagram is commuting
$$\xymatrix
{0\ar[r]& (N,\alpha_N)\ar[r]^{j}\ar@{=}[d]_{\Id}& (A,\alpha_A)\ar@<0.5ex>[r]^-\tau& (L,\alpha_L)\ar[r]& 0\\
0\ar[r]& (N,\alpha_N)\ar[r]^\chi& (P,\alpha_P)\ar@<0.5ex>[r]^-{\overline{\tau}}\ar@{-->}[u]_{\overline{\pi}}& (K,\alpha_K)\ar@<0.5ex>[l]^-{\sigma}\ar[r]\ar[u]_\pi& 0.
}$$
Hence, the following  commutative diagram hold:
$$\xymatrix
{0\ar[r]& (M,\alpha_M)\ar[r]^{i}& (K,\alpha_K)\ar@<0.5ex>[r]^-\pi\ar@{-->}[d]^{\overline{\pi}\circ \sigma}& (L,\alpha_L)\ar[r]\ar@{=}[d]^{\Id}& 0\\
0\ar[r]& (N,\alpha_N)\ar[r]& (A,\alpha_A)\ar@<0.5ex>[r]^-\tau & (L,\alpha_L)\ar[r]& 0,
}$$
i.e., there is a map $\overline{\pi}\circ \sigma$ such that $\pi=\tau\circ \overline{\pi} \circ \sigma$.
Lemma \ref{lem213} end the proof.

 c) and d) For a Hom-preLie algebra $(L, \alpha_L)$, we  consider the homology chain complex $CL_{\star}^{\alpha}(L)$, which is  $CL_{\star}^{\alpha}(L, \mathbb{K})$ and $\mathbb{K}$ is endowed with the trivial  Hom-co-representation structure.

As $\mathbb{K}$-vector spaces, let $I_L$ be the  subspace of $L \otimes L$ spanned by the elements of the form $\alpha_L(x_1)\otimes x_2x_3- x_1x_2\otimes \alpha_L(x_3)-\alpha_L(x_2)\otimes x_1x_3+ x_2x_1\otimes \alpha_L(x_3) , x_1, x_2, x_3 \in L$. That is, $I_L = {\rm Im}\ \left( d_3 : CL_3^{\alpha}(L) \to CL_2^{\alpha}(L) \right)$.

Now we denote the quotient  $\mathbb{K}$-vector space   $L \otimes L/ I_L$  by $\frak{uce}(L)$. The class $x_1 \otimes x_2 + I_L$ is denoted  $\{x_1,x_2\}$, for $x_1, x_2 \in L$.

One can see from the above discussion, the following identity holds:
\begin{equation}
\{\alpha_L(x_1),x_2x_3\} - \{x_1x_2, \alpha_L(x_3)\} - \{\alpha_L(x_2),x_1x_3\}+\{x_2x_1,\alpha_L(x_3)\}=0,
\end{equation}
for all $x_1, x_2, x_3 \in L$.

Now $d_2(I_L)=0$, so  $d_2$ induces a $\mathbb{K}$-linear map $u_L : \frak{uce}(L) \to L$, given by $u_L(\{x_1,x_2\})=x_1x_2$. Moreover $(\frak{uce}(L), \widetilde{\alpha})$, where $\widetilde{\alpha} : \frak{uce}(L) \to \frak{uce}(L)$ is defined by $\widetilde{\alpha}(\{x_1,x_2\}) = \{\alpha_L(x_1), \alpha_L(x_2) \}$,  is a  Hom-preLie algebra with respect to the bracket $\{x_1,x_2\}\{y_1,y_2\}= \{x_1x_2,y_1y_2\}$ and $u_L : (\frak{uce}(L), \widetilde{\alpha}) \to (L,\alpha_L)$ is a homomorphism of Hom-preLie algebras. Actually, $\rm{Im}$$u_L = LL$, and $(L,\alpha_L)$ is a perfect Hom-preLie algebra, so $u_L$ is a surjective homomorphism.

From the construction, it follows that  ${\Ker}(u_L) = HL_2^{\alpha}(L)$, so we have the extension
$$0 \to (HL_2^{\alpha}(L), \widetilde{\alpha}_{\mid}) \to (\frak{uce}(L), \widetilde{\alpha}) \stackrel{u_L}\to (L,\alpha_L) \to 0$$
which is central, since ${\rm Ker}(u_L)\frak{uce}(L) = \frak{uce}(L){\rm Ker}(u_L)=0$. It is universal, since for any central extension $0 \to (M,\alpha_M) \to (K,\alpha_K) \stackrel{\pi} \to (L,\alpha_L) \to 0$  there exists the homomorphism of Hom-preLie algebras  $\varphi : (\frak{uce}(L), \widetilde{\alpha}) \to (K,\alpha_K)$ given by $\varphi(\{x_1,x_2\})=k_1k_2, \pi(k_i)=x_i, i = 1, 2$, such that $\pi \circ \varphi = u_L$. Moreover, $(\frak{uce}(L), \widetilde{\alpha})$ is a perfect Hom-preLie algebra, so by Lemma \ref{lem213}, $\varphi$ is unique.
\epf

\bcor\label{cor211}
a) Let $0 \to (M, \alpha_M) \stackrel{i} \to (K,\alpha_K) \stackrel{\pi} \to (L, \alpha_L) \to 0$ be a universal   $\alpha$-central extension, then $HL_1^{\alpha}(K) = HL_2^{\alpha}(K) = 0$.

 b)  Let $0 \to (M, \alpha_M) \stackrel{i} \to (K,\alpha_K) \stackrel{\pi} \to (L, \alpha_L) \to 0$  be a central extension such that $HL_1^{\alpha}(K) = HL_2^{\alpha}(K) = 0$, then $0 \to (M, \alpha_M) \stackrel{i} \to (K,\alpha_K) \stackrel{\pi} \to (L, \alpha_L) \to 0$  is a universal  central extension.
\ecor
\bpf
{\it a)} If $0 \to (M, \alpha_M) \stackrel{i} \to (K,\alpha_K) \stackrel{\pi} \to (L, \alpha_L) \to 0$ is a universal  $\alpha$-central extension, then $(K,\alpha_K)$ is perfect by Remark \ref{re212} and Lemma \ref{lem214}, so  $HL_1^{\alpha}(K) = 0$. By Lemma \ref{lem215} and Theorem \ref{thm211} {\it c)} and {\it d)} the universal central extension corresponding to $(K,\alpha_K)$ is split, so $HL_2^{\alpha}(K) = 0$.

{\it b)} $HL_1^{\alpha}(K) = 0$ implies that $(K,\alpha_K)$ is a perfect Hom-preLie algebra.
$HL_2^{\alpha}(K) = 0$ implies that $(\frak{uce}(K),\widetilde{\alpha}) \stackrel{\sim} \to (K,\alpha_K)$. So there exists a map $\sigma:(K,\alpha_K) \to (\frak{uce}(K),\tilde{\alpha})$ such that  $u_K \circ \sigma={\Id}$.
Consider the central extension
$0 \to (P,\alpha_P) \to(N,\alpha_N) \stackrel{t} \to (K,\alpha_K)\to 0.$
Since the universal of the extension
$0 \to(HL_\alpha^2(K),\widetilde{\alpha_\mid}) \to(\frak{uce}(K),\tilde{\alpha}) \stackrel{u_K} \to(K,\alpha_K)\to 0$, there exists a unique homomorphism $\varphi:(\frak{uce}(K),\tilde{\alpha})\to (N,\alpha_N)$ such that the following diagram
$$\xymatrix
{0\ar[r]& (HL_\alpha^2(K),\widetilde{\alpha_\mid})\ar[r]& (\frak{uce}(K),\tilde{\alpha})\ar@<0.5ex>[r]^-{u_K}\ar@{-->}[d]^{\exists 1 \varphi}& (K,\alpha_K)\ar@<0.5ex>[l]^-\sigma\ar[r]\ar@{=}[d]^{\Id}& 0\\
0\ar[r]& (P,\alpha_P)\ar[r]& (N,\alpha_N)\ar@<0.5ex>[r]^-t & (K,\alpha_K)\ar[r]& 0
}$$
is commutative,
i.e., $t\circ \varphi \circ \sigma=u_K \circ \sigma={\Id}$. So every  central extension of $(K,\alpha_K)$ is split.
 Theorem \ref{thm211} {\it  b)} ends the proof.
\epf

\bdefn\rm\label{def211}
An $\alpha$-central extension $0 \to (M, \alpha_M) \stackrel{i} \to (K,\alpha_K) \stackrel{\pi} \to (L, \alpha_L) \to 0$ is said to be universal if for any central extension $0 \to (R, \alpha_R) \stackrel{j} \to (A,\alpha_A) \stackrel{\tau} \to (L, \alpha_L) \to 0$  there exists a unique homomorphism $\varphi : (K,\alpha_K) \to (A, \alpha_A)$ such that $\tau \circ \varphi = \pi$.

An $\alpha$-central extension $0 \to (M, \alpha_M) \stackrel{i} \to (K,\alpha_K) \stackrel{\pi} \to (L, \alpha_L) \to 0$ is said to be $\alpha$-universal if for any $\alpha$-central extension $0 \to (R, \alpha_R) \stackrel{j} \to (A,\alpha_A) \stackrel{\tau} \to (L, \alpha_L) \to 0$  there exists a unique homomorphism $\psi : (K,\alpha_K) \to (A, \alpha_A)$ such that $\tau \circ \psi = \pi$.
\edefn

\bre\rm
 Obviously, every $\alpha$-universal $\alpha$-central extension is an $\alpha$-central
extension which is universal in the sense of Definition \ref{def211}. In case $\alpha_M = \Id$
both notions coincide with the definition of universal central extension.
\ere

\bprop\label{prop213}
a) Let $0 \to (M, \alpha_M) \stackrel{i} \to (K,\alpha_K) \stackrel{\pi} \to (L, \alpha_L) \to 0$  and  $0 \to (N, \alpha_N) \stackrel{j} \to (F,\alpha_F) \stackrel{\rho} \to (K, \alpha_K) \to 0$  be central extensions.  If $0 \to (N, \alpha_N) \stackrel{j} \to (F,\alpha_F) \stackrel{\rho} \to (K, \alpha_K) \to 0$ is a universal central extension, then $0 \to (P, \alpha_P) = {\rm Ker} (\pi \circ \rho) \stackrel{\chi} \to (F,\alpha_F) \stackrel{\pi \circ \rho} \to (L, \alpha_L) \to 0$ is an  $\alpha$-central extension which is universal in the sense of Definition \ref{def211}.

b) Let $0 \to (M, \alpha_M) \stackrel{i} \to (K,\alpha_K) \stackrel{\pi} \to (L, \alpha_L) \to 0$  and  $0 \to (N, \alpha_N) \stackrel{j} \to (F,\alpha_F) \stackrel{\rho} \to (K, \alpha_K) \to 0$  be central extensions with $(F,\alpha_F)$ a perfect Hom-preLie algebra. If $0 \to (P, \alpha_P) = {\rm Ker} (\pi \circ \rho) \stackrel{\chi} \to (F,\alpha_F) \stackrel{\pi \circ \rho} \to (L, \alpha_L) \to 0$ is an $\alpha$-universal $\alpha$-central extension, then $0 \to (N, \alpha_N) \stackrel{j} \to (F,\alpha_F) \stackrel{\rho} \to (K, \alpha_K) \to 0$  is a universal central extension.
\eprop
\bpf
 a) If $0 \to (N, \alpha_N) \stackrel{j} \to (F,\alpha_F) \stackrel{\rho} \to (K, \alpha_K) \to 0$ is a universal central extension, then $(F,\alpha_F)$ and $(K,\alpha_K)$ are perfect Hom-preLie algebras by Lemma \ref{lem214}. Moreover,
$0 \to (P, \alpha_P) = {\rm Ker} (\pi \circ \rho) \stackrel{\chi} \to (F,\alpha_F) \stackrel{\pi \circ \rho} \to (L, \alpha_L) \to 0$ is an $\alpha$-central extension by Lemma \ref{lem215}.

 In order to obtain the universality, for any central extension $0 \to (R, \alpha_R) \to (A,\alpha_A) \stackrel{\tau} \to (L, \alpha_L) \to 0$ construct the pull-back extension  $0 \to (R, \alpha_R) \to (K\times_L A,\alpha_K\times \alpha_A) \stackrel{\bar{\tau}} \to (K, \alpha_K) \to 0$, where $K \times A=\{(k,a) \in K \times A \mid \tau(a)=\pi(k) \}$ satisfies the following diagram commutate
 $$\xymatrix
{0\ar[r]& (R,\alpha_R)\ar[r]^{j}\ar@{=}[d]_{\Id}& (A,\alpha_A)\ar@<0.5ex>[r]^-\tau& (L,\alpha_L)\ar[r]& 0\\
0\ar[r]& (R,\alpha_R)\ar[r]^-\chi& (K\times_L A,\alpha_K\times \alpha_A)\ar@<0.5ex>[r]^-{\overline{\tau}}\ar@{-->}[u]_{\overline{\pi}}& (K,\alpha_K)\ar[r]\ar[u]_\pi& 0
}$$
  Since $0 \to (N, \alpha_N) \stackrel{j} \to (F,\alpha_F) \stackrel{\rho} \to (K, \alpha_K) \to 0$ is a universal central extension, there exists a unique
homomorphism $\varphi:(F,\alpha_F)\to (K\times_L A,\alpha_K\times \alpha_A)$ such that $\tau\circ \bar{\pi} \circ \varphi=\pi\circ\rho$. Then the homomorphism $\bar{\pi} \circ \varphi$ satisfies that $\tau\circ \bar{\pi} \circ\varphi=\pi \circ \rho$ and it is unique by Lemma \ref{lem213}.

b) The fact that  $(F,\alpha_F)$ is perfect implies that $(K,\alpha_K)$ is perfect by Lemma \ref{lem212}. In order to prove
the universality of the central extension  $0 \to (N, \alpha_N) \stackrel{j} \to (F,\alpha_F) \stackrel{\rho} \to (K, \alpha_K) \to 0$, let us consider any central extension $0 \to (R, \alpha_R) \to (A,\alpha_A) \stackrel{\sigma} \to (K, \alpha_K) \to 0$, then $0 \to {\rm Ker} (\pi \circ \sigma)  \to (A,\alpha_A) \stackrel{\pi \circ \sigma} \to (L, \alpha_L) \to 0$ is an $\alpha$-central extension by
Lemma \ref{lem215}.

The $\alpha$-universality of $0 \to (P, \alpha_P) = {\rm Ker} (\pi \circ \rho) \stackrel{\chi} \to (F,\alpha_F) \stackrel{\pi \circ \rho} \to (L, \alpha_L) \to 0$ implies the existence of a unique homomorphism $\omega :(F,\alpha_F)\to (A,\alpha_A)$ such that $\pi\circ \sigma \circ \omega=\pi\circ \rho$.

Since $(F,\alpha_F)$ is a perfect Hom-preLie algebra and ${\rm{Im}}(\sigma\circ \omega)\subseteq {\rm{Im}}(\rho)+{\Ker}(\pi)$, we have $\sigma\circ \omega=\rho$.
Lemma \ref{lem213} end the proof.
\epf

\bdefn\rm
A Hom-preLie algebra $(L,\alpha_L)$ is said to be $\alpha$-perfect if
$$L=\alpha_L(L)\alpha_L(L).$$
\edefn

\bre\rm
a) When $\alpha_L = \Id$, the notion of $\alpha$-perfect coincides with the notion of perfection.

b) Obviously, if $(L,\alpha_L)$ is an $\alpha$-perfect Hom-preLie algebra, then it is perfect. Nevertheless the converse is not true in
general. For example, the three-dimensional Hom-preLie algebra $(L,\alpha_L)$ with the basis ${a_1,a_2,a_3}$, the multiplication given by $a_1a_2=-a_2a_1=a_3, a_1a_3=-a_3a_1=a_2, a_2a_3=-a_3a_2=a_1$ and the endomorphism $\alpha_L=0$ is perfect, but it is not $\alpha$-perfect.

c) If $(L,\alpha_L)$ is $\alpha$-perfect, then $L = \alpha_L(L)$, i.e., $\alpha_L$ is surjective. Nevertheless the converse is not true. For instance, the two-dimensional Hom-preLie algebra with the basis $\{a_1,a_2\}$, the multiplication given by $a_1a_1=a_1,a_1a_2=a_2a_1=a_2a_2=0$  and the endomorphism $\alpha_L$ represented by the matrix
$\big(\begin{smallmatrix} 1 & 0 \\
0 & 2 \end{smallmatrix}  \big).$ Obviously the
homomorphism $\alpha_L$ is surjective, but $\alpha_L(L)\alpha_L(L)=\langle a_1 \rangle.$
\ere

\blem\label{lem216}
 Let $0 \to (M, \alpha_M) \stackrel{i} \to (K,\alpha_K) \stackrel{\pi} \to (L, \alpha_L) \to 0$ be a central extension and $(K,\alpha_K)$ an $\alpha$-perfect Hom-preLie algebra. If there exists a
homomorphism of Hom-preLie algebras $f : (K, \alpha_K) \to (A, \alpha_A)$ such that $\tau \circ f=\pi$, where $0 \to (N, \alpha_N) \stackrel{j} \to (A,\alpha_A) \stackrel{\tau} \to (L, \alpha_L) \to 0$ is an $\alpha$-central extension, then $f$ is unique.
\elem
\bpf
The proof is similar to Lemma \ref{lem213}.
\epf

\bthm\label{thm212}
An $\alpha$-perfect Hom-preLie algebra admits a universal $\alpha$-central extension.
\ethm
\bpf
For an $\alpha$-perfect Hom-Leibniz algebra $(L, \alpha_L)$, we consider the chain complex
$CL_\star^\alpha(L, \K)$, where $\K$ is considered with trivial Hom-co-representation structure. Take the quotient vector space $\frak{uce}_\alpha(L)=\frac{\alpha_L(L)\otimes \alpha_L(L)}{I_L}$, where $I_L$ is the vector
subspace of $\alpha_L(L)\otimes \alpha_L(L)$ generated by the elements of the form
$$\alpha_L(x_1)\otimes x_2x_3- x_1x_2\otimes \alpha_L(x_3)-\alpha_L(x_2)\otimes x_1x_3+ x_2x_1\otimes \alpha_L(x_3).$$

Observe that every summand of the form $\alpha_L(x_1)\otimes x_2x_3$ or $x_1x_2\otimes \alpha_L(x_3)$
belongs to $\alpha_L(L)\otimes \alpha_L(L)$, since $L$ is $\alpha$-perfect,  $x_1x_2\in L = \alpha_L(L)\alpha_L(L)\subseteq \alpha_L(L)$.

We denote by $\{\alpha_L(x_1), \alpha_L(x_2)\}$ the equivalence class $\alpha_L(x_1)\otimes \alpha_L(x_2)+I_L$. $\frak{uce}_\alpha(L)$ is endowed with a structure of Hom-preLie algebra with respect to
the multiplication
$$\{\alpha_L(x_1), \alpha_L(x_2)\}\{\alpha_L(y_1), \alpha_L(y_2)\}=\{\alpha_L(x_1)\alpha_L(x_2),\alpha_L(y_1)\alpha_L(y_1)\}$$
and the homomorphism $\overline{\alpha}:\frak{uce}_\alpha(L) \to \frak{uce}_\alpha(L)$ defined by
$$\overline{\alpha}(\{\alpha_L(x_1), \alpha_L(x_2)\})=\{\alpha^2_L(x_1), \alpha^2_L(x_2)\}.$$

The restriction of the differential $d_2 : CL^\alpha_2 (L) = L\otimes L \to CL^\alpha_1 (L) = L$ to
$\alpha_L(L)\otimes \alpha_L(L)$ vanishes on $I_L$, so it induces a linear map $U_\alpha : \frak{uce}_\alpha(L) \to L$,
that is given by $U_\alpha(\{\alpha_L(x_1), \alpha_L(x_2)\}) = \alpha_L(x_1)\alpha_L(x_2)$. Moreover, thanks to be
$(L, \alpha_L)$ $\alpha$-perfect,  $U_\alpha$ is a surjective homomorphism, since ${\rm{Im}}(U_\alpha)=\alpha_L(L)\alpha_L(L)=L,$
\begin{align*}
&U_\alpha(\{\alpha_L(x_1), \alpha_L(x_2)\}\{\alpha_L(y_1), \alpha_L(y_2)\})=U_\alpha\{\alpha_L(x_1)\alpha_L(x_2),\alpha_L(y_1)\alpha_L(y_1)\}\\
=&(\alpha_L(x_1)\alpha_L(x_2))(\alpha_L(y_1)\alpha_L(y_1))=U_\alpha\{\alpha_L(x_1), \alpha_L(x_2)\}U_\alpha\{\alpha_L(y_1), \alpha_L(y_2)\}
\end{align*}
and
\begin{align*}
&\alpha_L\circ U_\alpha\{\alpha_L(x_1),\alpha_L(x_2)\}=\alpha_L(\alpha_L(x_1)\alpha_L(x_2))=\alpha^2_L(x_1)\alpha^2_L(x_2)\\
=&U_\alpha\{\alpha^2_L(x_1),\alpha^2_L(x_2)\}=U_\alpha\circ \overline{\alpha}\{\alpha_L(x_1),\alpha_L(x_2)\}.
\end{align*}
So we have constructed the extension
$$0 \to ({\Ker}(U_\alpha), \overline{\alpha}_{\mid}) \to (\frak{uce}_{\alpha}(L), \overline{\alpha}) \stackrel{U_\alpha}\to (L,\alpha_L) \to 0,$$
which is central. Indeed, for any $\{\alpha_L(x_1),\alpha_L(x_2)\}\in {\Ker}(U_\alpha), \{\alpha_L(x_1),\alpha_L(x_2)\}\in \frak{uce}_{\alpha}(L)$, we have that
$$\{\alpha_L(x_1), \alpha_L(x_2)\}\{\alpha_L(y_1), \alpha_L(y_2)\}=\{\alpha_L(x_1)\alpha_L(x_2),\alpha_L(y_1)\alpha_L(y_1)\}=0$$
since $\{\alpha_L(x_1), \alpha_L(x_2)\}\in {\Ker}(U_\alpha) \Leftrightarrow \alpha_L(x_1)\alpha_L(x_2)=0$.

 In a similar way it is verified that $\frak{uce}_{\alpha}(L){\Ker}(U_\alpha) = 0$.

Finally, this central extension is universal $\alpha$-central. Indeed, consider an $\alpha$-central extension $0 \to (M, \alpha_M) \stackrel{j} \to (K,\alpha_K) \stackrel{\pi} \to (L, \alpha_L) \to 0$. We define the map $\Phi:\frak{uce}_{\alpha}(L) \to (K,\alpha_K)$ by $\Phi(\{\alpha_L(x_1),\alpha_L(x_2)\})=\alpha_K(k_1)\alpha_K(k_2)$, where $\pi(k_i)=x_i,x_i\in L,k_i\in K,i=1,2.$ Now we are going to check that $\Phi$ is a
homomorphism of Hom-preLie algebras such that $\pi\circ \Phi=U_\alpha$. Indeed:
\begin{enumerate}
\item[1)] $\Phi$ is well defined: if $\pi(k_i)=x_i=\pi(k'_i),i=1,2$, then $k_i-k'_i\in {\Ker}(\pi)=M$, so $\alpha_K(k_i)-\alpha_K(k'_i)\in \alpha_M(M)$. Hence,
\begin{align*}
    &\Phi(\{\alpha_L(x_1),\alpha_L(x_2)\})=\alpha_K(k_1)\alpha_K(k_2)=(\alpha_K(k'_1)+\alpha_M(m_1))(\alpha_K(k'_2)+\alpha_M(m_2))\\
    =&\alpha_K(k'_1)\alpha_K(k'_2)+\alpha_K(k'_1)\alpha_M(m_2)+\alpha_M(m_1)\alpha_K(k'_2)+\alpha_M(m_1)\alpha_M(m_2)\\
    =&\alpha_K(k'_1)\alpha_K(k'_2).
\end{align*}
\item[2)] $\Phi$ keeps multiplication:
\begin{align*}
   &\Phi(\{\alpha_L(x_1),\alpha_L(x_2)\}\{\alpha_L(y_1),\alpha_L(y_2)\})=\Phi\{\alpha_L(x_1)\alpha_L(x_2),\alpha_L(y_1)\alpha_L(y_2)\}\\
     =&(\alpha_K(k_1)\alpha_K(k_2))(\alpha_L(l_1)\alpha_L(l_2))
    =\Phi\{\alpha_L(x_1),\alpha_L(x_2)\}\Phi\{\alpha_L(y_1),\alpha_L(y_2)\},
\end{align*}
where $\pi(k_i)=x_i,\pi(l_i)=y_i,i=1,2.$
\item[3)] $\alpha_K\circ \Phi=\Phi\circ \overline{\alpha}$:
\begin{align*}
    \alpha_K\circ\Phi(\{\alpha_L(x_1),\alpha_L(x_2)\})=&\alpha_K(\alpha_K(k_1)\alpha_K(k_2))=\alpha^2_K(k_1)\alpha^2_K(k_2)\\
    =&\Phi\{\alpha^2_L(x_1),\alpha^2_L(x_2)\}=\Phi\circ \overline{\alpha}\{\alpha_L(x_1),\alpha_L(x_2)\}.
\end{align*}
\item[4)] $\pi\circ \Phi=U_\alpha$:
\begin{align*}
&\pi\circ \Phi \{\alpha_L(x_1),\alpha_L(x_2)\}=\pi(\alpha_K(k_1)\alpha_K(k_2))=(\pi\circ\alpha_K(k_1))( \pi\circ\alpha_K(k_2))\\
=&(\alpha_L\circ\pi(k_1))(\alpha_L\circ\pi(k_2))=\alpha_L(x_1)\alpha_L(x_2)=U_\alpha\{\alpha_L(x_1),\alpha_L(x_2)\}.
\end{align*}
\end{enumerate}

To end the proof, we must verify the uniqueness of $\Phi$. First at all, we check that
$\frak{uce}_{\alpha}(L)$ is $\alpha$-perfect. Indeed,
$$\overline{\alpha}(\{\alpha_L(x_1),\alpha_L(x_2)\})\overline{\alpha}(\{\alpha_L(y_1),\alpha_L(y_2)\})=\{\alpha^2_L(x_1)\alpha^2_L(x_2),\alpha^2_L(y_1)\alpha^2_L(y_2)\},$$
i.e., $\overline{\alpha}(\frak{uce}_{\alpha}(L))\overline{\alpha}(\frak{uce}_{\alpha}(L))\subseteq \frak{uce}_{\alpha}(L).$

For the converse inclusion, by $\alpha_L(x_i)\in L = \alpha_L(L)\alpha_L(L)$,
we have that:
\begin{align*}
    \{\alpha_L(x_1),\alpha_L(x_2)\}=&\left\{\alpha_L\left(\sum_i\lambda_i(\alpha_L(l_{i_1})\alpha_L(l_{i_2}))\right),\alpha_L\left(\sum_j\mu_j(\alpha_L(l'_{j_1})\alpha_L(l'_{j_2}))\right)\right\}\\
    =&\sum\lambda_i\mu_j\{\alpha^2_L(l_{i_1})\alpha^2_L(l_{i_2}),\alpha^2_L(l'_{j_1})\alpha^2_L(l'_{j_2})\}\\
    =&\sum\lambda_i\mu_j(\{\alpha^2_L(l_{i_1}),\alpha^2_L(l_{i_2})\}\{\alpha^2_L(l'_{j_1}),\alpha^2_L(l'_{j_2})\})\\
    =&\sum\lambda_i\mu_j\overline{\alpha}\{\alpha_L(l_{i_1}),\alpha_L(l_{i_2})\}\overline{\alpha}\{\alpha_L'(l_{j_1}),\alpha_L'(l_{j_2})\}
    \in \overline{\alpha}(\frak{uce}_{\alpha}(L))\overline{\alpha}(\frak{uce}_{\alpha}(L)).
\end{align*}
Now Lemma \ref{lem216} ends the proof.
\epf


\begin{thebibliography}{99}
\bibitem{Casas2} J. Casas, M. Ladra, Stem extensions and stem covers of Leibniz algebras. Georgian Math. J. 9 (2002), 659--669.
\bibitem{Casas4} J. Casas, A. Vieites, Central extensions of perfect of Leibniz algebras. In: Recent Advances in Lie Theory: Research and Exposition in Mathematics 25 (2002), 189--196.
\bibitem{Casas3} J. Casas, M. Ladra, Computing low dimensional Leibniz homology of some
perfect Leibniz algebras. Southeast Asian Bulletin of Mathematics 31 (2007), 683--690.
\bibitem{Casas1} J. Casas, N. Corral, On Universal Central Extensions of Leibniz Algebras. Comm. Algebra 37 (2009), 2104--2120.
\bibitem{Casas5}J. Casas, M. Insua,  N. Pacheco Rego, On universal central extensions of Hom-Leibniz algebras. J. Algebra Appl. 13 (2014), 1450053, 22pp.
\bibitem{Casas6} J. Casas, N. Pacheco Rego, On the universal $\alpha$-central extension of the
semi-direct product of Hom-Leibniz algebras. arXiv:1309.5204.
\bibitem{Khmaladze} X.
 Garc\'{\i}a-Mart\'{\i}nez,
 E. Khmaladze, M. Ladra,
Non-abelian tensor product and homology of Lie superalgebras.
J. Algebra 440 (2015), 464--488.

\bibitem{Gnedbaye1} A. Gnedbaye, Third homology groups of universal central extensions of a
Lie algebra. Afrika Math. 10 (1999), 46--63.
\bibitem{Gnedbaye2} A. Gnedbaye, A non-abelian tensor product of Leibniz algebras. Ann. Inst. Fourier 49 (1999), 1149--1177.

\bibitem{Kurdiani} R. Kurdiani, T. Pirashvili, A Leibniz algebras structure on the second tensor
power. J. Lie Theory 12 (2002), 583--596.
\bibitem{Makhlouf1} A. Makhlouf, S. Silvestrov, Hom-algebra structures. J. Gen. Lie Theory Appl. 2 (2008), 51--64.

\bibitem{Yau} D. Yau, Hom-Novikov algebras. J. Phys. A 44 (2011), 085202, 20 pp.





\end{thebibliography}
\end{document}